\documentclass[12pt]{amsart}  
\usepackage{graphicx,psfrag}
\usepackage{color}
\usepackage{array}
\usepackage{colortbl}
\usepackage{amsthm,amssymb,hyperref}
\usepackage{tabularx}
\newcommand{\ben}{\begin{enumerate}}
\newcommand{\een}{\end{enumerate}}
\newcommand{\ble}{\begin{lemma}}
\newcommand{\ele}{\end{lemma}}
\newcommand{\bth}{\begin{theorem}}
\renewcommand{\eth}{\end{theorem}}
\newcommand{\bpr}{\begin{prop}}
\newcommand{\epr}{\end{prop}}
\newcommand{\bco}{\begin{corollary}}
\newcommand{\eco}{\end{corollary}}
\newcommand{\bcon}{\begin{conjecture}}
\newcommand{\econ}{\end{conjecture}}
\newcommand{\bex}{\begin{exa}}
\newcommand{\eex}{\end{exa}}
\newcommand{\barr}{\begin{array}}
\newcommand{\earr}{\end{array}}
\newcommand{\btab}{\begin{tabular}}
\newcommand{\etab}{\end{tabular}}
\newcommand{\beq}{\begin{equation}}
\newcommand{\eeq}{\end{equation}}
\newcommand{\bea}{\begin{eqnarray*}}
\newcommand{\eea}{\end{eqnarray*}}
\newcommand{\bal}{\begin{align*}}
\newcommand{\bce}{\begin{center}}
\newcommand{\ece}{\end{center}}
\newcommand{\bpi}{\begin{picture}}
\newcommand{\epi}{\end{picture}}
\newcommand{\bpp}{\begin{picture}}
\newcommand{\epp}{\end{picture}}
\newcommand{\bfi}{\begin{figure} \begin{center}}
\newcommand{\efi}{\end{center} \end{figure}}
\newcommand{\bprf}{\begin{proof}}
\newcommand{\eprf}{\end{proof}\medskip}

\newcommand{\bsl}{\begin{slide}{}}
\newcommand{\esl}{\end{slide}}
\newcommand{\bfr}{\begin{frame}}
\newcommand{\efr}{\end{frame}}

\newcommand{\hqed}{\hfill \qed}

\newcommand{\ol}{\overline}

\newcommand{\hso}[1]{\hspace{-1pt}}

\newcommand{\qmq}[1]{\quad\mbox{#1}\quad}

\newcommand{\emp}{\emptyset}

\newcommand{\sbe}{\subseteq}

\newcommand{\case}[4]{\left\{\barr{ll}#1&\mbox{#2}\\#3&\mbox{#4}\earr\right.}
\newcommand{\fl}[1]{\lfloor #1 \rfloor}
\newcommand{\ce}[1]{\lceil #1 \rceil}

\def\<{\langle}
\def\>{\rangle}

\newcommand{\ree}[1]{(\ref{#1})}

\renewcommand{\th}{\theta}

\newcommand{\Si}{\Sigma}

\newcommand{\bbZ}{{\mathbb Z}}
\newcommand{\cA}{{\mathcal A}}

\newcommand{\cT}{{\mathcal T}}

\newcommand{\ab}{\ol{a}}
\newcommand{\Ab}{\ol{A}}

\newcommand{\Rb}{\ol{R}}
\newcommand{\Sb}{\ol{S}}

\DeclareMathOperator{\Av}{Av}

\DeclareMathOperator{\diag}{diag}

\newcommand{\sv}{Springer-Verlag Lecture Notes in Math.\/}

\newtheorem{theorem}{Theorem}
\newtheorem{corollary}[theorem]{Corollary}
\newtheorem{conjecture}[theorem]{Conjecture}

\newtheorem{lemma}[theorem]{Lemma}
\newtheorem{prop}[theorem]{Proposition}

\usepackage{amsmath}
\usepackage{amsfonts}
\newfont{\bb}{msbm10}

\def\:{\! :\!}

\newcommand\rb{\ol{r}}
\renewcommand\sb{\ol{s}}

\newcommand\dna{\downarrow}
\newcommand\mv{\check{m}}
\newcommand\nv{\check{n}}
\newcommand\rv{\check{r}}
\renewcommand\sv{\check{s}}
\newcommand\Rv{\check{R}}
\newcommand\Sv{\check{S}}
\renewcommand\Av{\check{A}}

\newcommand\mh{\hat{m}}
\newcommand\nh{\hat{n}}

\DeclareMathOperator{\antidiag}{antidiag}

\usepackage{tikz}

\begin{document}

\title{Dihedral Transportation and $(0,1)$-Matrix Classes}

 \author{Richard A. Brualdi
 \and
 Bruce Sagan
 }

\maketitle

 \begin{abstract}  Let $R$ and $S$ be two vectors of real numbers whose entries have the same sum.  In the transportation problems one wishes to find a matrix $A$ with row sum vector $R$ and column sum vector $S$.  If, in addition, the two vectors only contain nonnegative integers then one wants the same to be true for $A$.  This can always be done and the transportation algorithm gives a method for explicitly calculating $A$.  We can restrict things even further and insist that $A$ have only entries zero and one.  In this case, the Gale-Ryser Theorem gives necessary and sufficient conditions for $A$ to exist and this result can be proved constructively.  One can let the dihedral group $D_4$ of the square act on matrices.  Then a subgroup of $D_4$ defines a set of matrices invariant under the subgroup.  So one can consider analogues of the transportation and $(0,1)$ problems for these sets of matrices.  For every subgroup, we give conditions equivalent to the existence of the desired type of matrix.


\end{abstract}

\section{Introduction}

Let $D_4$ be the dihedral group of the square.  Write $\rho_\th$ for rotation counter-clockwise through $\th$ radians and $r_m$ for reflection in a line of slope $m$.  Then 
$$
D_4=\{\rho_0,\rho_{\pi/2},\rho_\pi,\rho_{3\pi/2},r_0,r_{+1},r_{-1},r_\infty\}.
$$
The non-identity elements of $D_4$ are uniquely identified by their subscripts, and we let $D_b\le D_4$ be the cyclic subgroup generated by the element with subscript $b$.  There are also two subgroups of $D_4$ isomorphic to the Klein $4$-group, namely
$$
D_\times =\{\rho_0,\rho_\pi,r_{+1},r_{-1}\}
$$
and
$$
D_+ = \{\rho_0,\rho_\pi,r_0,r_\infty\}.
$$
The subscripts of $D_\times$ and $D_+$ are mnemonic, geometrically representing the two reflection lines in each subgroup.  A complete list of non-identity subgroups of $D_4$ is
$$
 D_{\pi/2}=D_{3\pi/2},  D_\pi, D_0, D_{+1}, D_{-1}, D_\infty, D_\times, D_+, D_4.
$$

For each of these subgroups $D_b$ (now including $D_\times$, $D_+$, and $D_4$) acting on $m\times n$ matrices (where it is implicitly assumed that $m=n$ if one of $\rho_{\pi/2}$, $r_{+1}$ or $r_{-1}$ is in $D_b$), we consider the transportation (both real and integral) and $(0,1)$-problems for those matrices invariant under $D_b$. We call the resulting classes of matrices {\it dihedral matrix classes}. The cases $D_{\pi}$ and  $D_\times$ were considered in a paper of Brualdi and Ma~\cite{RABMA}.  The invariant matrices for $D_{\pi}$ are  the so-called  {\it centrosymmetric} matrices. Since $D_{\pi}$ is a subgroup of $D_\times$, the invariant matrices for $D_\times$ are also centrosymmetric. As pointed out in \cite{RABMA}, there are centrosymmetric matrices that are not invariant under $D_\times$. For example,  the matrix
\[\left[\begin{array}{cccc}
0&0&1&0\\
1&0&0&0\\
0&0&0&1\\
0&1&0&0\end{array}\right]\]
is centrosymmetric but  is not invariant under either of the two reflections $r_{+1}$ and  $r_{-1}$.

Given a real matrix $A$ we let $R=R(A)$ and $S=S(A)$ be the row sum and column sum vectors of $A$ with components $r_i=r_i(A)$ and $s_j=s_j(A)$, respectively.
We let $\cT(R,S)$ denote the corresponding {\em transportation class} which consists  of all nonnegative real matrices with row sum vector $R$ and column sum vector $S$. We  also use the notation
$$
\cT^b(R,S) = \{A\in\cT(R,S)\ |\ D_b A = A\}
$$
and
$$
\cT^b_\bbZ(R,S) = \{A\in\cT^b(R,S)\ |\  A\in\bbZ^{m\times n}\}.
$$
For the $(0,1)$-problem,  ${\mathcal A}(R,S)$ and ${\mathcal A^b}(R,S)$ denote the subsets of $\cT(R,S)$ and $\cT^b(R,S)$, respectively, whose entries are 0 and 1.
In all cases we assume, without specific mention, the obvious necessary condition for our classes to be nonempty, namely that $\Si R =\Si S$  where, for any matrix $X$, $\Si X$ is the sum of the entries of $X$.  We  assume, also without specific mention, that in discussing $\cT^b_\bbZ(R,S)$ and  $\cA^b(R,S)$, the vectors $R,S$ have nonnegative integral components.  Finally, for $\cA^b(R,S)$, we always assume that $R$ and $S$ have no component bigger than $n$ and $m$, respectively.

Recall that we can obtain an element $T\in \cT(R,S)$ by letting
\beq
\label{tij}
t_{i,j} = \frac{r_i s_j}{N}
\eeq
where $N=\Sigma R = \Sigma S$.  

If we wish to construct a matrix $T\in\cT_\bbZ(R,S)$, then we can use the transportation algorithm.  Pick any $r_i$ and $s_j$.  If $r_i\le s_j$ then let $t_{i,j}=r_i$, remove the $i$th row of $T$ and the corresponding component of $R$, and replace $S$ by the vector obtained by decreasing its $j$th component by $s_j$.  If $s_j\le r_i$ then we apply the same construction with the roles of the rows and columns reversed.  (If $r_i=s_j$ it does not matter which of the two possibilities we use.)  We then iterate the process until all row and column sums are as they should be.

For $\cA(R,S)$ one must be more careful.
Given a nonnegative integral vector $R$, we let $R^\dna$ denote the weakly decreasing rearrangement of $R$, and we let $R^*$ denote the conjugate of $R^\dna$ viewed as an integer partition.  Note that $R^*$ is weakly decreasing by definition.  Given two weakly decreasing  vectors $R=(r_1,r_2,\ldots,r_m)$ and $S=(s_1,s_2,\ldots,s_n)$, we say $R$ {\em majorizes} $S$ and write $R\succeq S$, if for all indices $\ell$
\beq
\label{maj}
r_1+r_2+\dots+r_\ell\ge s_1+s_2+\dots+s_\ell
\eeq
and $\Si R = \Si S$.  We  also write $S\preceq R$ and say that $S$ is {\it majorized} by $R$. If $R,S$ are not necessarily weakly decreasing, then we define $R\succeq S$ (or $S\preceq R$) to mean $R^\dna\succeq S^\dna$.
The Gale-Ryser theorem (see e.g. \cite{RAB}) asserts that $\cA(R,S)\ne\emptyset$ if and only 
\beq
\label{GR}
S\preceq R^*\quad \mbox{(the {\it Gale-Ryser condition}).}
\eeq
If (\ref{GR}) holds,
then we can construct an element $A\in\cA(R,S)$ using the Gale-Ryser algorithm as follows.
\ben
\item  Pick any $j$ and set the entries in column $j$ with the largest $s_j$ row sums equal to one and the rest of the entries equal to zero, breaking ties arbitarily,
\item Replace $R$ by the vector obtained by decreasing its  largest $s_j$ entries by one (using  tie breaking as determined in (1)).  Replace $S$ by the vector obtained by removing $s_j$ and return to the first step until both vectors are zeroed out.
\een

It  will  be helpful to have the following notation.  For a nonegative integer $n$, let
$$
\nv = \fl{n/2}\qmq{and}\nh=\ce{n/2}.
$$
Also, if $A$ is a matrix, then $R_i$ and $S_j$ will always denote the $i$th row and $j$th column of $A$, respectively.

Our goal in this paper is to determine under what conditions the various dihedral matrix classes, as determined by the subgroups of $D_4$, are nonempty.


\section{The rotation $\rho_\pi$}

As mentioned in the introduction, these centrosymmetric matrices were considered in~\cite{RABMA}.  So here we content ourselves with stating their results. In order to state them more clearly, we assume some obvious necessary conditions.
Clearly a matrix invariant under $\rho_\pi$ must have palindromic row and column sum vectors.
We say that a palindromic vector $R=(r_1,r_2,\ldots,r_n)$ is {\em initially nonincreasing} provided that $r_1\ge r_2\ge\cdots\ge r_{\nv}$. By permuting within upper rows and within lower rows, and similarly for the columns, a centrosymmetric matrix can always be assumed to have initially-nonincreasing row and column sum vectors.

\bth
\label{Tpi}
We  have $\cT^{\pi}(R,S)\neq\emp$ if and only if $R$ and $S$ are palindromic.
The same is true for $\cT_\bbZ^{\pi}(R,S)$.\hqed
\eth

\bth
\label{Api}
\begin{itemize}
\item[(i)] Let $m$ and $n$ be even. Then $\cA^{\pi}(R,S)\neq\emp$ if and only if $R$ and $S$ are palindromic and $S\preceq R^*$.
\item[(ii)] Let $m$ be odd and $n$ be even, the case where $m$ is even and $n$ odd being similar.  Assume that $R$ and $S$ are initially nonincreasing, palindromic vectors with $r_{\mh}$ even. Let vectors $R'$ and $S'$ be obtained, respectively,  by deleting $r_{\mh}$ from $R$  and  by decreasing by one the first and last $r_{\mh}/2$ entries of $S$. 
Then $\cA^{\pi}(R,S)\neq\emp$ if and only if $\cA^{\pi}(R',S')\neq\emp$.
\item[(iii)] Let $m$ and $n$ both be odd, and assume that $R$ and $S$ are initially nonincreasing, palindromic vectors with $r_{\mh}$ and $s_{\nh}$ of the same parity. Let  vectors $R'$ and $S'$ be obtained, respectively,  by deleting $r_{\mh}$ and
by decreasing by $1$ the first and last $\lfloor s_{\nh}/2\rfloor$ entries of $R$,  and  by deleting $s_{\nh}$, and by decreasing by $1$ the first and last $\lfloor r_{\mh}/2\rfloor$
entries of $S$. Then $\cA^{\pi}(R,S)\neq\emp$ if and only if $\cA^{\pi}(R',S')\neq\emp$.\hqed
\end{itemize}
\eth


\section{The reflections $r_{-1}$ and $r_{+1}$}

In this section we will consider the subgroups $D_{-1}, D_{+1},$ and $D_{\times}$ generated by  the reflections $r_{-1}$ and/or $r_{+1}$. 

\bth
\label{T-1}
We  have $\cT^{-1}(R,S)\neq\emp$ if and only if $R=S$
The same is true for $\cT^{-1}_\bbZ(R,S)$.
\eth
\bprf
The proofs for the arbitrary and integral cases are the same.
To see the forward implication, it suffices to observe  that $r_{-1}$, which is ordinary matrix transposition,  interchanges the row and column sum vectors of a matrix.
For the reverse, merely note that if  $R=S$ then the diagonal matrix $\diag(r_1,\dots,r_n)$ provides a desired matrix.
\eprf

Given a vector $S=(s_1,s_2,\dots,s_n)$, we denote its reversal by
$$
S^r = (s_n,\dots,s_2,s_1).
$$
The next result follows from Theorem~\ref{T-1} and the fact that if $r_{+1}A=A$ if and only if $A$ can be obtained by rotation through $\pi/2$ radians of a matrix $A'$ with $r_{-1}A'=A'$ (i.e. transposition with respect  to the antidiagonal).
\bth
\label{T+1}
We  have $\cT^{+1}(R,S)\neq\emp$ if and only if $S=R^r$
The same
 is true for $\cT^{+1}_\bbZ(R,S)$.\hqed
\eth

Now we consider what happens for the subgroup $D_\times =\{\rho_0,\rho_\pi,r_{+1},r_{-1}\}$.
\bth
\label{Tx}
We have $\cT^\times(R,S)\neq\emp$ if and only if
	\ben
	\item[(a)]  $R=S$, and 
	\item[(b)] $R$ is palindromic.
	\een
The same is true for $\cT^\times_\bbZ(R,S)$.
\eth
\bprf
We will do both the arbitrary and integral cases at the same time.  The forward direction follows immediately from Theorems~\ref{T-1} and~\ref{T+1}.  On the other hand, if we are given (a) and (b) then it is easy to verify that 
\beq
\label{diag}
A=\diag(r_1/2,\dots,r_n/2)+\antidiag(r_1/2,\dots,r_n/2)
\eeq
is an element in $\cT^\times(R,S)$.  And for $\cT^\times_\bbZ(R,S)$ one merely rounds up the elements in the diagonal matrix  and rounds down those in the antidiagonal matrix.
\eprf

We now deal with the case of $(0,1)$-matrices.  For $r_{-1}$ this follows from a result of Fulkerson, Hoffman, and McAndrew~\cite{FHM} .  See~\cite[pp.\ 179--182]{RAB} for details.
\bth
\label{A-1}
We  have $\cA^{-1}(R,S)\neq\emp$ if and only if $R=S$ and $R\preceq R^*$.\hqed
\eth

Note that Theorem \ref{A-1} is equivalent to the fact  that, for $R=S$, there is a symmetric matrix in $\cA(R,R)$ if and only if $\cA(R,R)\ne\emptyset$.

The following result follows from the previous one in the same way that Theorem~\ref{T+1} follows from Theorem~\ref{T-1}.
\bth
\label{A+1}
We  have $\cA^{+1}(R,S)\neq\emp$ if and only if $S=R^r$ and $R^r\preceq R^*$.\hqed
\eth

The nonemptiness of $\cA^{\times}(R,R)$ was characterized in \cite{RABMA} as follows. 

\bth\label{new}
We have $\cA^{\times}(R,R)\ne\emptyset$ if and only if $\cA^{\pi}(R,R)\neq\emp$. \hqed
\eth

Recall that the characterization for $\cA^\pi(R,S)$, and thus for $\cA^{\times}(R,R)$ is given in Theorem~\ref{Api}.


\section{The reflections $r_\infty$ and  $r_0$}

In this section we will consider the subgroups generated by  the reflections $r_\infty$ and/or $r_0$.  First, however, we introduce some useful notation.
Call an integral  matrix $A$ {\em even} if all its entries are even.   Also let $o(A)$ be the number of odd entries of $A$.
Given an integral vector $R$ and an odd positive integer $n$, we define $A^R$ to be the $m\times n$ (0,1)-matrix whose only nonzero entries are $a^R_{i,\nh}$ for the indices $i$ such that $r_i$ is odd.  Given an integral vector $S$ and odd positive integer $m$, we define $A^S$ in a similar way.   Finally given $R,S$ and  both $m$ and $n$ are odd we define $A^+$  by
\beq
\label{a+}
a^+_{i,j} = \max\{a^R_{i,j},a^S_{i,j}\}.
\eeq
In other words, $A^+=A^R+A^S$ except in the case when the central elements of both $R$ and $S$ are odd in which case the central entry of the sum is too large by one.

\bth
\label{Tinfty}
\ben
\item[(I)] We have $\cT^\infty(R,S)\neq\emp$ if and only if
	\ben
	\item[(a)]  $S$ is palindromic.
	\een
\item[(II)] We have $\cT^\infty_{\bbZ}(R,S)\neq\emp$ if and only if  {\rm (a)} is true  and 
	\ben
	\item[(b)] if $n$ is even then $R$ is even, and if $n$ is odd then  $s_{\nh}\ge o(R)$. 
	\een
\een
\eth
\bprf
(I) For the forward implication, take $A$ such that $r_\infty A = A$.  Since $r_\infty$ exchanges columns equidistant from the vertical mid-line of $A$, we must have that $S$ is a palindromic.  For the other direction, it suffices to show that equation~\ree{tij} defines a matrix with palindromic $S$-vector.  Indeed, using the fact that $S$ is palindromic, 
$$
t_{i,m-j+1}= \frac{r_i s_{m-j+1}}{N}=\frac{r_i s_j}{N} = t_{i,j}.
$$

\medskip

(II)    First we note that  if  $r_\infty A = A$ then $a_{i,j}=a_{i,n-j+1}$ for all $i,j$.  Thus when $n$ is even every element in the $i$th row is repeated twice and $R$ is even.  On the other hand, if $n$ is odd then $r_i$ is odd if and only if $a_{i,\nh}$ is odd.  This gives  the inequality in (b).

For the reverse implication, we modify the transportation matrix algorithm as follows.  Let $\Rb=R-R(A^R)$ and $\Sb=S-S(A^R)$.
Note that $\Rb$ is even by definition of $A^R$ and $\Sb$ still has nonnegative entries because of (b).
Construct $\Ab\in \cT^\infty_{\bbZ}(\Rb,\Sb)$ by  letting $\ab_{1,1}=\ab_{1,n}=\min\{\rb_1/2,\sb_1\}$ and applying recursion.
Now form $A\in \cT^\infty_{\bbZ}(R,S)$ by adding one to the $\ab_{i,\nh}$ for all $i$ such that $r_i$ is odd.
\eprf

The next result follows from Theorem~\ref{Tinfty} in the same way that Theorem~\ref{T+1} follows from Theorem~\ref{T-1}.
\bth
\label{T0}
\ben
\item[(I)] We have $\cT^0(R,S)\neq\emp$ if and only if
	\ben
	\item[(a)]  $R$ is palindromic.
	\een
\item[(II)] We have $\cT^0_{\bbZ}(R,S)\neq\emp$ if and only if {\rm (a)} is true and 
	\ben
	\item[(b)] if $m$ is even then $S$ is even, and  if $m$ is odd then $r_{\mh}\ge o(S)$. \hqed 
	\een
\een
\eth

We now consider the subgroup $D^+=\{\rho_0,\rho_\pi,r_0,r_\infty\}$.
\bth
\label{T+}
We have $\cT^+(R,S)\neq\emp$ if and only if $\cT^\infty(R,S)\neq\emp$ and $\cT^0(R,S)\neq\emp$.
The same is true in the integral case.
\eth
\bprf
The forward directions follow immediately from the fact that $\cT^+(R,S)=\cT^\infty(R,S)\cap\cT^0(R,S)$.  The converse for $\cT^+(R,S)$ is proved in the usual way using~\ree{tij}.  For $\cT^+_\bbZ(R,S)$, we use a method similar to the one  given in the proof of Theorem~\ref{Tinfty}.  We consider the vectors $\Rb=R-R(A^+)$ and $\Sb=S-S(A^+)$.  We then construct a matrix $\Ab$ by making assignments $\ab_{1,1}=\ab_{1,n}=\ab_{m,1}=\ab_{m,n}=\min\{\rb_1/2,\sb_1/2\}$ and recursing.  Finally, we let $A=\Ab+A^+$.  
\eprf

\bth
\label{Ainfty}
We have $\cA^\infty(R,S)\neq \emp$ if and only if  conditions {\rm (a)} and {\rm (b)} from Theorem~$\ref{Tinfty}$ are satisfied as well as
\ben
\item[(c)] $\Sb\preceq \Rb^*$ where $\Rb$ is obtained from $R$ by subtracting one from every odd component and $\Sb$ is $S$ if $n$ is even or  $S$ with column $S_{\nh}$ removed if $n$ is odd.
\een
\eth
\bprf
Clearly if $A\in\cA^\infty(R,S)$ then it must satisfy the  two conditions from Theorem~\ref{Tinfty}.    
If $n$ is even then $\Rb=R$ and $\Sb=S$ so that  $\Rb^*\succeq \Sb$ by the Gale-Ryser Theorem.
In $n$ is odd, note that the ones in column $S_{\nh}$ must occur exactly in the rows with odd sums.  Removing this column, we obtain a matrix $\Ab$ with $\Rb$ and $\Sb$ as its row and column vector.  Since such a matrix exists, we must have $\Rb^*\succeq \Sb$ by the Gale-Ryser Theorem again.

For the converse we have two cases.  First suppose that $n$ is even.  Then since $R$ is even we must have every element of $R^*$ repeated twice.  Let $R^*_1$  be $R^*$ where we only take one out of every pair of repeated elements.  Similarly, let 
$S_1=(s_1,\dots,s_{\nv})$.  Since $\Rb=R$ and $\Sb=S$, (c) implies that $S\preceq R^*$. It follows that $S_1\preceq R^*_1$.  Now use the Gale-Ryser algorithm to create a matrix $B\in\cA(R_1,S_1)$.  It follows that we have a  block matrix $A=[B\  r_\infty B]\in\cA^\infty(R,S)$.

Now consider the case when $n$ is odd.  Since $n-1$ is even, $\Rb$ is an even vector, $\Sb$ is palindromic, and $\Sb\preceq\Rb^*$ we can proceed as in the previous case to construct a matrix $\Ab\in\cA^\infty(\Rb,\Sb)$. Finally, we get the desired matrix $A$ by inserting  a middle column $S_{\nh}$ in $\Ab$ which has ones in exactly the rows of $R$ with odd sum.
\eprf

One might ask if (d) could be replaced by the ordinary Gale-Ryser condition $S\preceq R^*$.  But this condition is not strong enough to imply $\cA^\infty(R,S)\neq\emp$.  For an example of this, consider $R=(6,6,6,2,1,1)$ and $S=(4,4,2,2,2,4,4)$.  Clearly $S$ is palindromic and it is easy to check that $S\preceq R^*$.  Now suppose, towards a contradiction, that there exists $A\in\cA^\infty(R,S)$.  Form the matrix $\Ab$ as in the first paragraph of the preceding proof.  Then $\Ab$ has row and column vectors $\Rb=(6,6,6,2)$ and $\Sb=(4,4,2,2,4,4)$.  But $\Rb^*$ does not majorize $\Sb$ which contradicts the Gale-Ryser Theorem.

As with previous cases, the result for symmetry under $r_\infty$ is similar to the one for $r_0$.

\bth
\label{A0}
We have $\cA^0\neq \emp$ if and only if  conditions {\rm (a)} and {\rm (b)} from Theorem~$\ref{T0}$ are satisfied as well as
\ben
\item[(c)] $\Sb\preceq \Rb^*$ where $\Sb$ is obtained from $S$ by subtracting one from every odd component and $\Rb$ is $R$ if $m$ is even or  $R$ with column $R_{\nh}$ removed if $n$ is odd.\hqed
\een
\eth

Finally, we consider the $(0,1)$-case for $D_+$.

\bth
\label{A+}
We have $\cA^+(R,S)\neq \emp$ if and only if  conditions {\rm (a)} and {\rm (b)} from both Theorems~$\ref{Tinfty}$  and~$\ref{T0}$ are satisfied as well as
\ben
\item[(c)] if $n$ is odd then $o(R)=s_{\nh}$, if $m$ is odd then $o(S)=r_{\mh}$, and 
\item[(d)]  $\Sv\preceq \Rv^*$ where $\Rv=(\rv_1,\rv_2,\dots,\rv_{\mv})$  and $\Sv=(\sv_1,\sv_2,\dots,\sv_{\nv})$.
\een
\eth
\bprf
Suppose first that $A\in \cA^+(R,S)$.  Then clearly conditions {\rm (a)} and {\rm (b)} from both Theorems~\ref{Tinfty}  and~\ref{T0} are satisfied. 
To obtain (c) of the present result, note that condition (c) of Theorem~\ref{Ainfty} must also hold.  So, in particular, $\Si\Rb^*=\Si \Sb$ and this gives the desired equality when $n$ is odd.  The case when $m$ is odd follows similarly from Theorem~\ref{A0}.
 Finally, $\Rv$ and $\Sv$  are the row- and column-sum vectors for the submatrix  $\Av$ of $A$ sitting in the first $\mv$ rows and the first $\nv$ columns.  Thus  $\Rv^*\succeq \Sv$ follows from the Gale-Ryser Theorem.

For the converse, assume first that $m$ and $n$ are odd.  By condition (d) and the Gale-Ryser Theorem, we can construct an $\nv\times\mv$ matrix $\Av$ with row sum vector  $\Rv$ and column sum vector  $\Sv$.  Now the current condition (c) and condition (a) from Theorems~\ref{Tinfty}  and~\ref{T0} imply that there is an $\mv\times1$ matrix $B$, a $1\times\nv$ matrix $C$, and $a_{\mh,\nh}\in\{0,1\}$ such that the block matrix
$$
A=
\left[
\barr{ccc}
\Av 		&B			&r_\infty\Av\\
C 		&a_{\mh,\nh}	&r_\infty C\\
r_0 \Av	&r_0 B 		&\rho_\pi \Av
\earr
\right]
$$
is in $\cA^+(R,S)$.  If either $m$ or $n$ is even then condition (b) from Theorems~\ref{Tinfty}  and~\ref{T0} implies that deleting the appropriate row or column in $A$ above will give a matrix with the correct row and column sums to be in $\cA^+(R,S)$.
\eprf

\section{The case $D_{\pi/2}$}

We start, as usual, with the transportation problem.

\bth
\label{Tpi/2}
\ben
\item[(I)] We have $\cT^{\pi/2}(R,S)\neq\emp$ if and only if
	\ben
	\item[(a)]  $R=S$, and
	\item[(b)]  $R$ is palindromic.
	\een
\item[(II)] We have $\cT^{\pi/2}_{\bbZ}(R,S)\neq\emp$ if and only if  $R,S$ satisfy {\rm (a)} and {\rm (b)} as well as one of
	\ben
	\item[(c)] $r_1+r_2+\dots+r_{\nv}$ is even, or
	\item[(d)] $n$ is odd and $r_{\nh}\ge2$.
	\een
\een
\eth
\bprf
(I)  For the forward direction, suppose $A\in\cT^{\pi/2}(R,S)$.  Then  $\rho_{\pi/2} R_i=C_i$ which implies $R=S$.  And  $\rho_{\pi/2}^2 R_i = \rho_\pi R_i$ is $R_{n-i}$ read backwards so that (b) holds.

For the converse, it suffices to show that when (a) and (b) hold then the matrix defined by~\ree{tij} is invariant under $\rho_{\pi/2}$.  But this follows since
$$
t_{n-j+1,i} = \frac{r_{n-j+1} s_i}{N} = \frac{s_{n-j+1} r_i}{N} =\frac{r_i s_j}{N} = t_{i,j}.
$$

\medskip

(II)  We will first consider the case when $n$ is even. Given $A\in\cT^{\pi/2}_Z(R,S)$, we can write $A$ in the block form
\beq
\label{AB}
A= \left[ 
\barr{cc}
B & \rho_{\pi/2}^3 B\\
\rho_{\pi/2} B & \rho_{\pi/2}^2 B
\earr
\right]
\eeq
where $B$ is $\nv\times\nv$.  Since $R$ is palindromic by (a), it follows that
$$
r_1+r_2+\dots+r_{\nv}=\Sigma B + \Sigma (\rho_{\pi/2}^3 B) = 2\Sigma B
$$
so that (c) holds.

Now suppose, for $n$ still even, that we are given (a)--(c).    For any matrix $B$, the matrix $A=A(B)$ defined by~\ree{AB} is invariant under $\rho_{\pi/2}$.  Thus it suffices to show that we can define $B$ so that $A$ has the given row and column sums.  We will define $B=D+P$ where $D$ is a diagonal matrix and $P$ is a $(0,1)$-matrix with at most one $1$ in every row and column.  Define $D$ by $d_{i,i}=\rv_i$ for $1\le i\le \nv$.  It follows that $A(D)$ has rows sums $2\rv_i =r_i$ if $r_i$ is even or $r_i-1$ if $r_i$ is odd.  We use the matrix $P$ to correct for the odd row sums as follows.  Because of (c), there are an even number of $r_i$ which are odd,  $1\le i\le \nv$.  Let those $r_i$ be $r_{i_1},r_{i_2},\dots,r_{i_{2k}}$.  Let  $P$ be the $(0,1)$-matrix with $1$'s in positions $(i_1,i_2),\dots,(i_{2k-1},i_{2k})$.  Now $A=A(B)$ will have one added to row $i_{2j-1}$ by $B$ and to row $i_{2j}$ by $\rho_{\pi/2}^3 B$  for $1\le j\le k$ and similarly for the rows below the midpoint.  It follows that $A$ has the correct row sums and we are done with the case $n$ even.

We now deal with $n$ odd.  If $A\in\cT^{\pi/2}_{\bbZ}(R,S)$ then, similarly to the $n$ even case, we write
\beq
\label{ABC}
A= \left[ 
\barr{ccc}
B &C&\rho_{\pi/2}^3 B\\
\rho_{\pi/2} C&a_{\nh,\nh}&\rho_{\pi/2}^3 C\\
\rho_{\pi/2} B & \rho_{\pi/2}^2 C&\rho_{\pi/2}^2 B
\earr
\right]
\eeq
where $B$ is $\nv\times\nv$ and $C$ is  $\nv\times 1$.  If (c) holds, then we are done.  If not, then consider
$$
r_1+r_2+\dots+r_{\nv}= 2\Sigma B+\Sigma C.
$$
By our assumption about the left-hand side we must have $\Sigma C$ odd and so, in particular, $\Sigma C\ge1$.  But then
$$
r_{\nh} = 2\Sigma C + a_{\nh,\nh}\ge2
$$
and so (d) holds.

Finally, we must prove the converse when $n$ is odd.  If (c) holds, then we can construct the matrix $B$ as when $n$ is even, take $C$ to be a zero matrix, and set 
$a_{\nh,\nh}=r_{\nh}$ to obtain a matrix with the desired row and column sums.  If, instead, (d) holds then there are an odd number of $r_i$ which are odd,  $1\le i\le \nv$.  Let those $r_i$ be $r_{i_1},r_{i_2},\dots,r_{i_{2k+1}}$.  Construct that matrix $B$ as for $n$ even using $r_{i_1},r_{i_2},\dots,r_{i_{2k}}$.  Let $C$ be the matrix which is all zeros except for its $i_{2k+1}$ entry which is one.  And define $a_{\nh,\nh}=r_{\nh}-2\ge0$ by the assumption in (d).  It is now an easy matter to verify that we again have the desired sums in rows and columns.
\end{proof}

For the $(0,1)$-case we will need the following result of Brualdi and Ryser~\cite[Theorem 6.3.2]{RAB} about symmetric matrices whose entries are zeros, ones, and twos.
\bth
\label{B-R}
Let $R=(r_1,\dots,r_n)$ be a  vector of nonnegative integers.  There exists a symmetric $(0,1,2)$-matrix $M$ with row sum vector $R$ if and only if 
\beq
\label{sym}
2|I||J|\ge \sum_{i\in I} r_i - \sum_{j\not\in J} r_j
\eeq
for all $I,J\sbe\{1,2,\dots,n\}$.\hqed
\eth

We note that if in the previous theorem we have $R$ weakly decreasing (and the row vector of any symmetric matrix can be brought to this form by row and column interchanges), then it suffices to check the considerably smaller set of inequalities
$$
2kl\ge \sum_{i\le k} r_i -\sum_{i>l} r_i
$$ 
for all $1\le k\le j\le n$.

\bth
 We have $\cA^{\pi/2}(R,S)\neq\emp$ if and only if conditions (a)--(d) of Theorem~\ref{Tpi/2} hold and  $\Rb$ satisfies the inequalities~\ree{sym} where
$$
\Rb=\case{(r_1,r_2,\dots,r_{\nv})}{if $n$ is even,}{(r_1-1,r_2-1,\dots,r_s-1,r_{s+1},r_{s+2},\dots,r_{\nv})}{if $n$ is odd,}
$$
and $s=\fl{r_{\nh}/2}$.
\eth
\bprf
We begin with the case when $n$ is even.  Suppose first that $A\in\cA^{\pi/2}(R,S)$.  We have already shown that conditions (a)--(c) must be satisfied.  For the last condition, note that since $r_{\pi/2} A = A$ and $n$ is even this matrix must have the form~\ree{AB} for some $(0,1)$-matrix $B$.  It follows that $M=B+B^t$ is a symmetric $(0,1,2)$-matrix.  Furthermore, for $i\le n/2$ we have
\beq
\label{MA}
r_i(M)=r_i(B)+r_i(B^t)=r_i(B)+c_i(B)=r_i(B)+r_i(\rho_{\pi/2}^3B)=r_i(A).
\eeq
It follows from Theorem~\ref{B-R} that $\Rb$ must staisfy~\ree{sym}.

For the converse, using Theorem~\ref{B-R} again we may assume that there exists a symmetric $(0,1,2)$-matrix $M$ with $R(M)=\Rb$.  We claim that in fact there exists such an $M$ with no ones on the diagonal.  Indeed, using the symmetry of $M$ we have
$$
r_1+\dots+r_{n/2}=\Si M =2\sum_{i<j} m_{i,j} + \sum_i m_{i,i}.
$$
Since the left-hand side is even by condition (c), the same must be true of $\sum_i m_{i,i}$.  And because the only odd entries of $M$ are ones there must be an even number of them on $M$'s diagonal, say the entries $(i,i)$ for $i=i_1,i_2,\dots,i_{2k}$.  Consider the pair of ones on the diagonal in positions $i_{2j-1}$ and $i_{2j}$ for $1\le j\le k$.  Then there are three possibilities for the 
$2\times 2$ submatrix of $M$ in the rows and columns indexed by $i_{2j-1}$ and $i_{2j}$ depending on which of the three integers $0,1,2$ appear in the off-diagonal spots.  In each case, substitute the submatix on the left in the following table with the corresponding submatrix on the right.  It is easy to check that this does not change the row and column sums of $M$, and  now $M$ has only zeros and twos on the diagonal.

$$
\barr{c|c}
\text{initial submatrix}	&\text{substituted submatrix}\\
\hline
\left[\barr{cc}
1&0\\
0&1
\earr\right]
&
\left[\barr{cc}
0&1\\
1&0
\earr\right]
\rule{0pt}{30pt}
\\
\left[\barr{cc}
1&1\\
1&1
\earr\right]
&
\left[\barr{cc}
0&2\\
2&0
\earr\right]
\rule{0pt}{30pt}
\\
\left[\barr{cc}
1&2\\
2&1
\earr\right]
&
\left[\barr{cc}
2&1\\
1&2
\earr\right]
\rule{0pt}{30pt}
\earr
$$

We now write $M=B+B^t$ with the entries of the $(0,1)$-matrix $B$ defined as in the following chart for $i\le j$.  Note that from what we have just proved, if $m_{i,j}=m_{j,i}=1$ then we must actually have $i<j$.  

$$
\barr{c|c}
\text{entries of $M$}	&\text{entries of $B$}\\
\hline
m_{i,j}=m_{j,i}=0		&b_{i,j}=b_{j,i}=0\\
m_{i,j}=m_{j,i}=1		&b_{i,j}=0,\ b_{j,i}=1\\
m_{i,j}=m_{j,i}=2		&b_{i,j}=b_{j,i}=1
\earr
$$

Finally, we define $A$ using the matrix $B$ as in~\ree{AB}.  This matrix is clearly symmetric under $r_{\pi/2}$ and has the correct row and column sum vectors by conditions (a) and (b) and the equalities in~\ree{MA}.

Now suppose that $n$ is odd.  By interchanging rows and columns, we can assume that $R$ satisfies 
$r_1\ge r_2\ge\dots\ge r_{\nv}$.  Note that if there exists an $A\in\cA^{\pi/2}(R,S)$ then it must have the form given in~\ree{ABC}.  First  we claim that there is $A\in\cA^{\pi/2}(R,S)$ if and only if there is such a matrix where all the ones in $C$ precede all the zeros.  To prove the forward direction (the converse being trivial), suppose that the given matrix $A$ has a zero before a one in $C$. Without loss of generality we can assume the zero is in row $i$ and the one in row $i+1$.  But 
$r_i\ge r_{i+1}$ so that in some column of $A$ we must have a zero followed by a one  in these rows.  Suppose that this column is in $B$ as the case when it is in $\rho_{\pi/2}^3 B$ is similar.  So, taking account of symmetry, we have the situation depicted in (\ref{eq:newnew}) below:
\begin{equation}\label{eq:newnew}
A=\left[\begin{array}{cccccc|c|cccccc}
&&&&&&&&&&&&\\ 
&&&&1&&0&&&&0&1&\\ 
&&&&0&&1&&&&&&\\ 
&&&&&&&&&&&&\\ 
&&&&&&&&&&&&\\ 
&&&&&&&&&&&&\\ \hline
&0&1&&&&&&&&1&0&\\ \hline
&&&&&&&&&&&&\\ 
&1&0&&&&&&&&&&\\ 
&&&&&&&&&&&&\\ 
&&&&&&1&&&&&0&\\ 
&&&&&&0&&&&&1&\\ 
&&&&&&&&&&&&\end{array}\right].
\end{equation}
Now interchanging submatrices
$$
\left[\begin{array}{cc} 1 &0\\ 0&1\end{array}\right]
\leftrightarrow
\left[\begin{array}{cc} 0&1\\ 1&0\end{array}\right]
$$
in four different places maintains both the symmetry and the row sum vector while exchanging $a_{i,\nh}=0$ and $a_{i+1,\nh}=1$.  Continuing in this way we can put all the ones in $C$ before all the zeros.

Note that by its definition, $s$ is the number of ones in $C$.  So existence of  $A\in\cA^{\pi/2}(R,S)$ is equivalent to having such an $A$ with ones in the first $s$ rows of $C$ and zeros elsewhere in that submatrix.  Removing the central row and column of $A$, we see that this is equivalent to having a matrix with an even number of rows and columns which has $\Rb$ as the first half of its palindromic row sum vector, where $\Rb$ is as given in the statement of the theorem for $n$ odd.  So the case for $n$ odd reduces to the case when $n$ is even and we are done.
\eprf

\section{The group $D_4$ itself}

We finally deal with the full dihedral group.
\bth
\label{T4}  
\ben
\item[(I)]
We have $\cT^4(R,S)\neq\emp$ if and only if
	\ben
	\item[(a)]  $R=S$, and 
	\item[(b)] $R$ is palindromic.
	\een
\item[(II)]
We have $\cT^4_\bbZ(R,S)\neq\emp$ if and only if (a) and (b) hold as well as
	\ben
	\item[(c)] if $n$ is even then $R$ is even, and if $n$ is odd then  $r_{\nh}\ge o(R)$. 
	\een
\een
\eth
\bprf
(I)  The forward direction follows from Theorem~\ref{Tx} and the fact that $D_{\times}\subseteq D_4$.  For the reverse implication, it is easy to verify that if (a) and (b) are true then the matrix defined by~\ree{diag} is invariant under $D_4$.

\medskip

(II)  Similar to (I), the forward implication comes from Theorems~\ref{Tx} and~\ref{T+}.  For sufficiency, when $n$ is even we use~\ree{diag}.  When $n$ is odd, we let $\Av$ be the matrix defined as in~\ree{diag} but with all fractions rounded down.    It follows that $A=\Av+A^+$ is the desired matrix, where the entries of $A^+$ are defined by~\ree{a+}.
\eprf

For our final result, we characterized the  $(0,1)$-case.
\bth
\label{A4}
We have $\cA^4(R,S)\neq\emp$ if and only if conditions (a) and (b) of Theorem~\ref{T4} hold as well as
\ben
\item[(c)] if $n$ is even then $R$ is even, if $n$ is odd then $o(R)=r_{\nh}$, and 
\item[(d)]  $\Rv\preceq \Rv^*$ where $\Rv=(\rv_1,\rv_2,\dots,\rv_{\nv})$.
\een
\eth
\bprf
Necessity follows from the previous result and Theorem~\ref{A+}.  For the reverse implication, suppose first that $n$ is even.  By condition (d) and Theorem~\ref{A-1}, there is an $\nv\times\nv$ matrix $B$ with row and column sum vector $\Rv$ which is symmetric under matrix transposition.  It follows that the matrix $A$ defined by~\ref{AB} is invariant under $D_4$ and has the correct row and column sums by (c).  When $n$ is odd we construct $B$ as in the even case, then a matrix $\Av$ as in~\ref{ABC} where $C$ and $a_{\nh,\nh}$ are all zero, and finally let $A=\Av+A^+$ with entries given by~\ref{a+}.  Again, it is easy to see that $A$ has the desired properties.
\eprf

\nocite{*}
\bibliographystyle{alpha}

\bigskip

\begin{tabular}{ll}
\begin{minipage}[t]{8cm}
Richard A. Brualdi\\
Department of Mathematics\\
 University of Wisconsin\\
 Madison, WI 53706\\
 {\tt brualdi@math.wisc.edu}\\
 \end{minipage}
 &
 \begin{minipage}[t]{8cm}
 Bruce Sagan\\
 Department of Mathematics\\
 Michigan State University\\
 East Lansing, MI 48824\\
 {\tt sagan@math.msu.edu}
 \end{minipage}
 \end{tabular}

\end{document}